\documentclass{amsart}
\usepackage{amsgen}
\usepackage{amssymb}
\usepackage{amscd}
\newlength{\defbaselineskip} 

 \theoremstyle
 {plain}
\newtheorem{theo}{Theorem}[section]
 \newtheorem{lemm}[theo]{Lemma}
\newtheorem{propo}[theo]{Proposition}
\theoremstyle{definition}

\newtheorem{rema}[theo]{Remark}

 \numberwithin{equation}{section}
\numberwithin{equation}{section} \theoremstyle{definition}
\newtheorem{ex}{Example}[section] 
\begin{document}

\title{A cascade of determinantal Calabi--Yau threefolds}
\author{Grzegorz Kapustka}
\author{Micha\l\ Kapustka}
\maketitle
\begin{center} with Appendix by Piotr Pragacz
\end{center}
\begin{abstract} We study Kustin--Miller unprojections of Calabi--Yau threefolds. As an application we work out the geometric properties of
Calabi--Yau threefolds defined as linear sections of determinantal
 varieties.  We compute their Hodge numbers and describe the
 morphisms corresponding to the faces of their K\"{a}hler--Mori cone.
\end{abstract}
\section{Introduction}
By a Calabi--Yau threefold we mean a smooth complex projective
threefold with trivial canonical divisor, such that
$H^1(\mathcal{O}_X)=0$. The basic examples of such threefolds are
complete intersections in appropriate projective spaces or in
homogenous varieties. The aim of this paper is to enlarge the class
of easy to handle Calabi--Yau threefolds by working out the
geometric properties of Calabi--Yau threefolds whose ideals in
projective spaces are defined by the minors of an appropriate matrix
of linear forms. These varieties form a background for the
simultaneous application of various techniques concerning Calabi-Yau
threefolds, unprojections and the study of determinantal varieties
in general.

The considered determinantal Calabi--Yau threefolds were already
studied from a different point of view in \cite{GP} (see also
\cite{Bertin}). In particular, M.~Gross and S.~Popescu observed an
analogy between the descriptions of smooth del Pezzo surfaces $D'$
embedded by their anti-canonical divisors and descriptions of some
families of Calabi--Yau threefolds (we present it in Table
\ref{analogy}). By the symbols $X_{d_1,d_2,\dotsc}$, we denote
generic complete intersections of indicated degrees in the indicated
manifold.

 \begin{table}[h]\label{analogy}
 \caption{The analogy}
             \begin{center}
\vspace{-10pt}
\renewcommand*{\arraystretch}{1}
            \begin{tabular}{|c|c|c|}\hline
           $i$ & del Pezzo surfaces $D'$ & Calabi--Yau threefolds $X'$         \\ \hline

           1 &$D_6\subset \mathbb{P}(1,1,2,3)$& \\ \hline

           2 &$D_4\subset \mathbb{P}(1,1,1,2)$&$X_{8}\subset \mathbb{P}(1,1,1,1,4)$  \\ \hline

           3&$D_3\subset \mathbb{P}^3$& $X_5\subset \mathbb{P}^4$ \\ \hline

            4&$D_{2,2}\subset \mathbb{P}^4$& $X_{3,3}\subset \mathbb{P}^5$    \\ \hline
            5&$4 \times 4$ Pfaffians of a &$6 \times 6$ Pfaffians
            of a \\

             &$5\times 5$ skew-symmetric matrix & $7\times 7$ skew-symmetric matrix      \\ \hline
             6&$2\times 2$ minors of a $3\times 3$ matrix&$3\times 3$ minors of a $4\times 4$ matrix  \\ \hline
              7&$2\times 2$ minors of a $3\times4$ matrix&$3\times 3$ minors of a $4\times 5$ matrix  \\
              &obtained by deleting one row   &  obtained by deleting one row                \\
                  &from a  symmetric matrix    & from  a symmetric matrix                 \\ \hline

              8& $2\times 2$ minors of a $4\times4$ & $3\times 3$ minors of a $5\times
              5$ \\
                &symmetric matrix    &   symmetric matrix                 \\ \hline
 8'& $2\times 2$ minors of a $3\times 5$ &$3\times 3$ minors of a $4\times6$ \\
               &double-symmetric matrix  (see \ref{deg8})  &  double-symmetric matrix                  \\
               \hline\end{tabular}
\end{center}
          \end{table}

\vspace{0cm} We give a kind of explanation of this analogy by
joining the threefolds in Table \ref{analogy} by a natural
sequence (a "cascade") of conifold transitions corresponding to
projections (recall that Reid and Suzuki studied in \cite{Reid Su}
``cascades'' of del Pezzo surfaces). It is easy to observe that
lines in appropriate projective spaces can be described in the
following analogous way (see Table 2).
\begin{table}[h]\label{lines}
 \caption{lines}
             \begin{center}
\vspace{-10pt}
\renewcommand*{\arraystretch}{1}
            \begin{tabular}{|c|c|}\hline
           $i$ & lines      \\ \hline
           2 &$D_2\subset \mathbb{P}(1,1,2)$  \\ \hline

           3&$D_1\subset \mathbb{P}^2$ \\ \hline

            4&$D_{1,1}\subset \mathbb{P}^3$    \\ \hline
            5&$2 \times 2$ Pfaffians of a
            of a \\ & $3\times 3$ skew-symmetric matrix      \\ \hline
             6&$1\times 1$ minors of a $2\times 2$ matrix \\ \hline
              7&$1\times 1$ minors of a $2\times3$ matrix  \\
              &obtained by deleting one row                \\
                  &from a  symmetric matrix              \\ \hline

              8& $1\times 1$ minors of a $3\times3$  \\
                &symmetric matrix                  \\ \hline
 8'& $1\times 1$ minors of a $2\times 4$ \\
               &double-symmetric matrix      \\
               \hline\end{tabular}
\end{center}
\end{table}

We know that an anticanonically embedded del Pezzo surface of degree
$d$ is the projection of a del Pezzo surface of degree $d+1$ from a
generic point lying on the surface. Looking at these projections
from the opposite side we see a cascade of so-called Kustin--Miller
unprojections (see \cite{RP}) constructed as follows. We take a line
$l$ in some projective space and consider a generic del Pezzo
surface $D$ containing this line, then we find a meromorphic
function on $D$ with poles along $l$ and take the closure of its
graph. In this way we obtain a del Pezzo surface of degree higher by
one.

We can perform an analogous construction for del Pezzo surfaces
and Calabi-Yau threefolds (note that we can naturally generalize
it in higher dimensions). More precisely, we choose a del Pezzo
surface of degree $i$. We embed it into a Calabi--Yau threefold
$X'$ corresponding (in Table \ref{analogy}) to a del Pezzo surface
of degree $i-1$, and perform an unprojection. The difference
between the construction of a cascade of del Pezzo surfaces is
that both varieties (before the projection and after the
projection) are singular. We first need to prove that $X'$ may be
chosen to be nodal, which already appears to be technically
complicated (cf.~\cite[Thm.~2.1]{K}). In this part of the paper we
use the computer algebra system Singular and Lemma \ref{double
point st6} to handle specific examples. An alternative way is
however also shown in Theorem \ref{thm 4na4}. Next, by resolving
the singularities of $X'$ in such a way that the strict transform
$D'$ of $D$ is isomorphic to $D$, we obtain a Calabi--Yau
threefold $X$ containing a del Pezzo surface of degree $i$. We
next show that the surface $D$ can be contracted by an appropriate
linear system to a point that lies on a singular threefold. The
composition of the resolution and the contraction appears to be
the unprojection whose result belongs to the family of Calabi--Yau
threefolds corresponding to del Pezzo surfaces of degree $i$.

Using results of Namikawa \cite{namikawa1}, we obtain that the
Calabi--Yau threefold defined in $\mathbb{P}^7$ by the $3\times 3$
minors of a generic $4\times 4$ matrix of linear forms has Hodge
numbers $h^{1,1}=2$ and $h^{1,2}=34$. The Calabi--Yau threefolds
defined in $\mathbb{P}^8$ by the $3\times 3$ minors of a generic
$4\times 5$ partially symmetric matrix of linear forms in
$\mathbb{P}^8$ has Hodge numbers $h^{1,1}=2$ and $h^{1,2}=25$, and
the Calabi--Yau threefold defined in $\mathbb{P}^9$ by the
$3\times 3$ minors of a generic $5\times 5$ symmetric matrix of
linear forms in $\mathbb{P}^9$ has Hodge numbers $h^{1,1}=1$ and
$h^{1,2}=26$. Moreover, we give the descriptions of primitive
contractions corresponding to the faces of the K\"{a}hler--Mori
cones of these threefolds. We study the above varieties using the
Grassmann blow-up (see \cite{CM}) whose general properties are
discussed in Lemma \ref{michal}.  We also use results from
\cite{Co,JLP} and the Appendix by P.~Pragacz to this paper, all
concerning general determinantal varieties.

The results of this paper may help in studying applications of the
general theory of unprojections (see \cite{P,Reid 2}). In
particular we show a method of finding a nodal variety of given
type containing a Gorenstein variety. Such constructions were used
by Reid and Papadakis \cite{RP} to construct singular Fano
threefolds. In this paper we adapt the theory to the study of
Calabi-Yau threefolds. Our results suggest that it is natural to
try to connect all the families of Calabi--Yau threefolds using
unprojections (cf.~"Reid fantasy conjecture"). We hope to develop
these methods in order to construct and study many classes of
Calabi--Yau threefolds and unprojections between them including
unprojections of type II and higher. We also would like to use the
cascade to find mirrors to the considered determinantal
Calabi--Yau threefolds (cf.~\cite{B,Rod}).

\section*{Acknowledgements} We would like to express our gratitude to
P.~Pragacz for his encouragment and help. We thank L.~Borisov,
J.~Buczynski, S.~Cynk, A.~Langer for discussions, advice and
correcting mistakes.

\section{Del Pezzo of degree $\leq 5$}\label{1}
As a starting point let us consider the basic example of
Kustin--Miller unprojection (cf.~\cite{RP}). The del Pezzo surface
$D'\subset \mathbb{P}^3$ of degree $3$ is defined by the cubic
$c(x_1,x_2,x_3,x_4)=0$.
 We want to embed it into a singular Calabi--Yau threefold being a hypersurface
  $X'=X'_{8}\subset\mathbb{P}(1,1,1,1,4)$ (with coordinates $x_1,\dotsc,x_4,u$). Let hence $X'$ be defined by the equation
  $c(x_1,\dotsc,x_4)f(.)+ug(.)$, where $f$ an $g$ are a generic quintic and
a generic quartic in $\mathbb{P}(1,1,1,1,4)$. Clearly $X'$
contains the surface $D'$ in its natural embedding.
\begin{propo} The threefold $X'$ is a nodal Calabi--Yau threefold
with $60$ nodes. Moreover, the small resolutions of $X'$ have
Picard group of rank $2$.
\end{propo}
\begin{proof} This theorem can be proved by arguing as in Section 2.4 from
\cite{K}. Instead, let us show how the methods of this paper work
in this simplest example. First, the locus of singularities of
$X'$ is defined by the equations
$u=f(x_1,\dotsc,x_4,u)=g(x_1,\dotsc,x_4,u)=c(x_1,\dotsc,x_4)=0.$
After changing coordinates, the singularities are given locally by
$xy-zt=0$, hence these are $60$ ordinary double points. Let us now
describe the small resolution $X$ of $X'$. Consider the variety
$Y$ given in $\mathbb{P}(1,1,1,1,1,4)$ with coordinates
$(x_0,x_1,\dotsc,x_4,u)$ by the equations
\begin{center}$f(x_1,\dotsc,x_4,u)+x_0u$\\ and \\ $g(x_1,\dotsc,x_4,u)+x_0c(x_1,\dotsc,x_4)$.
\end{center}
The threefold $Y$ has only one singular point at $(1,0,0,0,0,0)$
being an isolated Gorenstein singularity. Indeed, since
$$u_1=g(x_1,\dotsc,x_4,u)+x_0c(x_1,\dotsc,x_4)$$ is a generic
quartic, we can change coordinates to $(x_0,\dotsc,x_4,u_1)$ and
see that $Y$ is a quintic in $\mathbb{P}^4$ with one singular
point, that has tangent cone given by $c(x_1,\dots,x_4)$.
Performing a weighted blow up $Y\subset \mathbb{P}(1,1,1,1,1,4)$
in the point $(1,0,\dotsc,0)$ with weights $1,1,1,1,4$, we obtain
its resolution of singularities $Z\subset \mathbb{P}(1,1,1,1,1,4)
\times\mathbb{P}(1,1,1,1,4)$. Observe that $Y$ contains $60$
lines, given by the equations
$$f(x_1,\dotsc,x_4,u)=g(x_1,\dotsc,x_4,u)=c(x_1,\dotsc,x_4)=u=0,$$
passing through the point $(1,0,0,0,0,0)\in Y$. These lines become
disjoint on $Z$. Now, $X$ is the image of the natural projection
$Z\rightarrow \mathbb{P}(1,1,1,1,4)$ and $Z\rightarrow X$ is a
small resolution that contracts exactly the considered $60$ lines.
 Finally, the rank of the
Picard group of $Y$ is $1$, thus the Picard group of $Z$, that is
the blow up of $Y$ in one point, has rank $2$ (see \cite[Prop.~2.6
II]{Har}).
\end{proof}

Using a different language, in the above proof we flopped the
exceptional curves of the blowing up of $D'\subset X'$, and
obtained a Calabi--Yau threefold containing a del Pezzo surface
$D$ of degree $3$.  The primitive contraction of type II with $D$
as exceptional divisor has image being a normal quintic in
$\mathbb{P}^4$. Observe that we found also in the proof the
quintic equation of the image.

Del Pezzo surfaces of degrees $4$ and $5$ have already been embedded into Calabi--Yau threefolds of
degrees $5$ and $9$ respectively in Section 5 from \cite{GM}. Observe however that using the method described
above we can find an exact description of the images also in these
cases. It remains to find appropriate embeddings for del Pezzo
surfaces of degree $\geq 6$.
\begin{rema} Observe that by analogy, an eventual Calabi--Yau threefold corresponding to a del
Pezzo surfaces of degrees $1$ would need to have degree $0$.
\end{rema}
\begin{rema} The descriptions from Table \ref{analogy} of the anti-canonical images of del Pezzo surfaces of degrees $i\leq 5 $ were
discussed in a more general context in \cite{GM1}. To see that in
the remaining cases the considered equations describe del Pezzo surfaces, it is
enough to show that these equations define smooth surfaces of
degree $i$ in $\mathbb{P}^i$ (see the Theorem of Nagata
\cite{nagata}).
\end{rema}
\section{Del Pezzo of degree 6}\label{st6}
Let $D'\subset \mathbb{P}^6$ be an anticanonically embedded del
Pezzo surface of degree 6. It can be proved that $D'$ is defined by
the $2 \times 2$ minors of a generic $3 \times 3$ matrix $M$ of
linear forms in $\mathbb{P}^6$.

We shall embed $D'$ into a nodal Calabi--Yau threefold defined by
the $6\times 6$ Pfaffians of a $7\times 7$ skew-symmetric matrix.
First, since each $n \times n$ matrix is a sum of a symmetric and
an skew-symmetric matrix, we can write $M=S+A$ as such a sum. Let
$B$ be the extra-symmetric (i.e. skew-symmetric and symmetric with
respect to the second diagonal) $6\times 6$ matrix
$$\left( \begin{array}{c}
\begin{array}{cc}
A&S^r\\
 -S^r&((A^r)^r)^t
\end{array}

\end{array}\right)$$

\vskip10pt where $.^r$ denotes the rotation of the matrix by $90$
degrees and $.^t$ is the transposition. Observe that the set of $2
\times 2$ minors of $M$ and the set
 of $4\times 4$ Pfaffians of $B$
 generate the same ideal (the ideal of $D'$). Let $C$ denote the $8 \times 8$ extra-symmetric matrix

$$\left( \begin{array}{c|c|c}
\begin{array}{c}
0
\end{array}
&
\begin{array}{ccc}
t_1&\dotsc&t_6
\end{array}
&
\begin{array}{c}
t_7
\end{array}

\\
\begin{array}{c}
-t_1\\
\vdots\\
-t_6
\end{array}
& \begin{array}{c}\text{\Huge{\emph{B}}}\end{array}&
\begin{array}{c}
t_6\\
\vdots\\
t_1
\end{array}\\
\begin{array}{c}
-t_7
\end{array}&\begin{array}{ccc}
-t_6&\dotsc&-t_1
\end{array}&0

\end{array}\right)$$

\vskip15pt where $t_1,..,t_7$ are linear forms in $\mathbb{P}^6$.
Let $C_1$ be the skew-symmetric matrix obtained from $C$ by deleting the
last row and the last column. Let $X'$ be a threefold defined by
$6\times 6$ Pfaffians of $C_1$ for a generic choice of $t_1,...,t_6$. To see that $X'$
contains $D'$, we use the fact that Pfaffians can be expanded
along any of their rows. In this way each $6\times 6 $ Pfaffian of
$C_1$ can be seen as an element of the ideal generated by the
$4\times 4$ Pfaffians of $B$.

Observe moreover that the $6\times 6$ Pfaffians of $C$ define a
smooth surface $G'$ of degree $20$ in $\mathbb{P}^6$, that is
contained in $X'$. Indeed, since $C$ can be represented in the
form
$$C=\left( \begin{array}{c}
\begin{array}{cc}
\tilde{A}&\tilde{S}^r\\
 -\tilde{S}^r&(((\tilde{A})^r)^r)^t
\end{array}

\end{array}\right)$$

where $\tilde{A}$ is a skew-symmetric and $\tilde{S}$ a symmetric
$4\times 4$ matrix, it follows that $G'$ is also defined by the
$3\times 3$ minors of the $4\times 4$ matrix $\tilde{A}+\tilde{S}$
(observe that in the upper left corner of this matrix we have
$t_7$). From the Porteous formula (see \cite{Fulton}) we obtain
that $G'$ has degree $20$. Moreover, since $\tilde{A}+\tilde{S}$
is a generic $4\times 4$ matrix with linear forms, $G'$ is smooth.

Now, using Singular \cite{GPS}, we compute (for some
$t_1,\dotsc,t_6$) that $G'$ and $D'$ have exactly $20$ points in
common (the radical of the ideal $\mathcal{I}_{G'}+
\mathcal{I}_{D'}$ has degree $20$). Moreover, the ideal
$\mathcal{I}_{G'}+ \mathcal{I}_{D'}$ has degree $20$. This means
that the ideal $\mathcal{I}_{G'}+ \mathcal{I}_{D'}$ is radical.
Hence we obtain two Weil divisor cutting transversally in each
singular point of $X'$. This suggests the following theorem.

\begin{theo}\label{nodal singular}  The threefold $X'\subset \mathbb{P}^6$
is a nodal Calabi--Yau threefold with $20$ nodes.
\end{theo}
\begin{proof}
Let $Y\subset \mathbb{P}^7$ (with coordinates $(x_0,\dotsc ,x_7))$
be a scheme defined by the vanishing of the $3\times3$ minors of
the $4\times 4$ matrix $F$ obtained from the matrix
$\tilde{A}+\tilde{S}$ (with linear forms depending of $7$
variables $x_1, \dotsc,x_7$) by replacing the form $t_7$ in its
upper left corner by the remaining free variable $x_0$. Notice
that the singular locus of the determinantal variety in
$\mathbb{P}^{15}$ parameterizing $4\times 4$ matrices of rank
$\leq 3$ is the locus of matrices of rank $\leq 2$. It follows now
from the Bertini theorem that $Y$ is smooth outside the point
$P=(1,0,\dotsc,0)$. By the discussion in \cite[p.~257]{Ha} we
conclude that the singularity of $Y$ has tangent cone being a cone
over a del Pezzo surface of degree $6$ determined by the vanishing
of the $2\times 2$ minors of the lower right $3\times 3$
sub-matrix of $F$. Moreover, from general properties of
determinantal varieties the singularity of $Y$ can be resolved by
blowing up $P$ (see Lemma \ref{michal}).

Let us consider the projection $\psi\colon \mathbb{P}^7 \supset Y-P
\rightarrow \mathbb{P}^6$ with center at the point $P$. The map
$\psi$ can be naturally extended to $\tilde{\psi}\colon
Y_P\rightarrow \mathbb{P}^6$, where $Y_P$ is the blow up of
$Y\subset \mathbb{P}^7$ in the point $P$. Denote by $D$ the
exceptional divisor of this blow-up. It is a del Pezzo surface of
degree $6$ (that is isomorphic to $D'$).
\begin{lemm}
 The morphism $\tilde{\psi}$ is birational, surjective onto $X'$ and injective outside the sum of $20$ lines
 contained in $Y_P$, which are contracted to $20$ points.
\end{lemm}
\begin{proof} Let us start with the proof that $\tilde{\psi}$ is birational.
To do this it is enough to prove that $\psi|_U$ is birational for
some open subset $U$ of $Y$. Choose a codimension $2$ linear space
$L\subset\mathbb{P}^7$ that does not contain $P$. Let $H$ be the
hyperplane spanned by $L$ and $P$. Let us take for $U$ above the
set $Y_L=Y\setminus H$. It is a subset of
$\mathbb{C}^7=\mathbb{C}^6\times \mathbb{C}$ with coordinates
chosen in such a way that the first $\mathbb{C}^6$ parameterizes
lines passing through the point $P$ and the last coordinate
parameterizes hyperplanes containing $L$. Recall moreover that
varying the upper left linear form $t_7$ in the matrix
$\tilde{A}+\tilde{S}$ we obtain a family $G'_{t_7}\subset X'$ of
surfaces. This permits us to find the following description of
$Y_L$ $$Y_L=\{ (p,x)\in \mathbb{C}^6\times \mathbb{C}\colon p\in
G'_{l_x}\},$$ where $l_x$ is the linear form defining the
hyperplane corresponding to $x$. In this notation $\psi|_{Y_L}$ is
given by the projection onto the $\mathbb{C}^6$. To study the
behavior of this projection consider the variety.
$$ \mathcal{Y}=\{ (p,l)\in \mathbb{P}^6\times \mathbb{C}^7\colon p\in G'_{l}\}.$$
We claim that for a fixed point $p\in X'- D'$, the set of linear
forms $t_7\in \mathbb{C}^7$ such that $p\in G_{t_7}$ is a
hyperplane. Indeed, the determinant of the matrix $B$ in the point
$p$ is zero. It follows that performing linear operations on rows
and columns of $B$ we can assume that the only nonzero entries in
$B(p)$ are contained in the $4\times 4$ sub-matrix obtained by
deleting the first and the last rows and columns. We obtain that
all $6\times 6$ Pfaffians of $C$ except one vanish in $p$. Thus,
this means that we have exactly one linear condition on the value
$t_7(x)$. Using the fact that $Y_L$ is a restriction of
$\mathcal{Y}$ obtained by choosing a line in $\mathbb{C}^7$ we
deduce that for a generic point $p$ of $X'\setminus D'$ there is
exactly one $x$ such that $G'_{l_x}\ni P$ (the point of
intersection of the hyperplane with the line). It follows that
$\tilde{\psi}$ is birational onto $X'$ and injective over
$X'\setminus D'$.

The same argument shows that the image $\tilde{\psi}(Y_P\setminus
D)$ is contained in the sum of all surfaces of the form $G'_l$. In
particular $$\tilde{\psi}(Y_P\setminus D)\subset(X'\setminus D')
\cup \bigcup_{l\in\mathbb{C}^7} D'\cap G'_l\,. $$ The set of all
intersection points of surfaces $G'_l$ with $D'$ is defined (set
theoretically) by the vanishing on $D'$ of the $3\times 3$ minors
of the matrix $T=\tilde{A}+\tilde{S}$ with $t_7=0$ . We continue by proving the following claim.\\
\emph{Claim} The $3\times 3$ minors of $T$ defines on $D'$ a set
of $20$ points.\\
Observe that if for a $3\times 3$ matrix with complex entries the
determinant is $0$, the upper left entry is $0$, the minor
obtained by deleting the first row and the first column is $0$ and
if the minor obtained by deleting the first row an the last column
is nonzero, then the minor obtained by deleting the last row an
the first column is zero. After some investigation this implies
that for a given $x$ the $3\times 3$ minors of $T$ and the
$2\times 2$ minors of the matrix obtained by deleting the first
row and the first column all vanish in $x$ if and only if either
all $2\times 2$ minors of the $3\times 4$ matrix obtained from $T$
by deleting the first row vanish in $x$, or all $2\times 2$ minors
of the $3\times 4$ matrix obtained from $T$ by deleting the first
column vanish in $x$. From the Porteous formula we conclude that
$T$ defines on $D'$ two disjoint sets of $10$ points. The claim
follows.

The above claim implies in particular that $\tilde{\psi}(D)=D'$.
Hence, since we know that both $D$ and $D'$ are del Pezzo surfaces
of degree $6$, we deduce that $\tilde{\psi}|_{D}$ is an
isomorphism onto $D'$. It remains to identify the set of curves
contracted by $\tilde{\psi}$. The only curves that might be
contracted by the considered morphism are curves that map to one
of the 20 intersection points of surfaces $G'_l$ with $D'$. An
easy calculation shows that the fibers over these points are one
dimensional (the fiber of $\mathcal{Y}\rightarrow \mathbb{P}^6$
under these points are $\mathbb{C}^7$).
\end{proof}

The above Lemma implies that $\tilde{\psi}$ is a small resolution
of $X'$, provided that $X'$ is normal.

\begin{lemm}
The variety $X'$ is normal.
\end{lemm}
\begin{proof}
Since $X'$ is of the expected codimension, we obtain by
general properties of Pfaffian subschemes that its singularities are
Gorenstein. We claim that the singularities of $X'$ are isolated.
Indeed, from the fact that the determinantal variety in
$\mathbb{P}^{14}$ determined by the $6\times 6$ Pfaffians of the
$7\times 7$ skew-symmetric matrix with the lower right $6\times 6$
sub-matrix being extra-symmetric has singular locus of
dimension $8$ and degree $20$ (calculation of the Jacobian ideal
with Singular \cite{GPS}) and since $X'$ is a linear section of this
variety, it follows from the Bertini theorem that the singularities
of $X'$ are isolated, Gorenstein, and in consequence normal.
\end{proof}
To conclude, observe that through each singular point of $X'$
there passes a family of smooth Weil divisors $G'_{t_7}$ cutting
transversally each other in this singular point. Indeed, it is
enough to prove that there exists a linear form $l$ such that
$G'_l$ and $D'$ meet transversally. We easily find such an $l$
using Singular or by calculation with a given example. Next, we
compute the partial derivatives of the expansions of the
determinants defining $G'_l$ along its first rows to see that
elements of the considered family meet transversely. The theorem
is now a direct consequence of the following Lemma.
\end{proof}
\begin{lemm}\label{double point st6}  Let $X'$ be a normal
threefold with only isolated Gorenstein threefold singularities
having a small resolution. Suppose that there exists seven Weil
divisors such that each two of them meet transversally in these
singular points, then $X'$ is a nodal threefold.
\end{lemm}
\begin{proof} Let us choose $P$ a singular point of $X'$. We shall first prove that the singularity of $X'$ in the point $P$ is of type $cA_1$
(i.e it can be described locally analytically by the equation
$x^2+y^2+z^2+tf(x,y,z,t)=0$). To do this take a generic hyperplane section $H$
through a singular point $P$ of $X'$. We need to prove that $P$ is
a singularity of type $A_1$ on this hyperplane section.

>From \cite[Cor.~1.12]{Reid 3} we know that the hyperplane section
of $P$ is a rational double point (an $ADE$ singularity). Suppose
that it is a singularity of type $A_k$, where $k\geq 2$. The
traces on $H$ of four of our Weil divisors $W_1, W_2, W_3$, and
$W_4$ are four curves cutting each other transversally. A
singularity of type $A_k$ can be described as the double covering
of $\mathbb{C}^2$ branched along the curve $x^2+y^k=0$. The image
of $W_1$, $W_2$, $W_3$, and $W_4$ by this covering are four curves
$C_1$, $C_2$, $C_3$, $C_4$ on $\mathbb{C}^2$ passing through
$(0,0)$. Observe that we can choose at least two of them (say
$C_1$ and $C_2$) not to be tangent to the line $x=0$. This follows
from the fact that $x=0$ is the only line passing through $(0,0)$
that splits in the double covering. We thus obtain a contradiction
since the Weil divisors mapping to $C_1$ and $C_2$ cannot cut
transversally. The case of singularities of type $D_n$ an $E_n$
can be treated similarly, and makes use of all seven divisors.

Since $P$ is a singularity of type $cA_1$ it can be described by
the equation $x^2+y^2+z^2+t^{2n}$. Consider its projective tangent
cone that is a quadric of rank $\geq 3$ that contain two disjoint
lines (corresponding to the tangent of our Weil divisors). It
follows that this quadric has rank $4$, so $P$ is an ordinary
double point.
\end{proof}
\begin{rema}
We can consider the morphism $\phi\colon X'\setminus D'\rightarrow
\mathbb{P}^7$ that associates to a point $x\in X'\setminus D'$ the
hyperplane of linear forms $l$ in $\mathbb{C}^7$ such that $x\in
G'_l$. Observe that $\phi$ is inverse to $\psi$.
\end{rema}
Let us see the above results in a different language. Blowing up
$D'$ we resolve the singularities of $X'$. Flopping the
exceptional divisors we obtain a smooth Calabi--Yau threefold $X$.
Denote by $D$ and $G$ the strict transforms to $X$ of $D'$ and
$G'$ respectively.
\begin{propo}\label{st6lemat} The linear
system $|G|$ defines a birational morphism $\pi\colon X\rightarrow
\tilde{X}$ with $D$ as exceptional locus into a normal variety
described by the vanishing of the $3\times 3$ minors of the
$4\times 4$ matrix $F$ (see the proof of Theorem \ref{nodal
singular}).
\end{propo}
\begin{proof} First, since each rational curve contracted by
$\tilde{\psi}\colon Y_P\rightarrow X'$ cuts $D''$ with
multiplicity $1$, one obtains that $\tilde{\psi}$ is exactly the
small resolution $X\rightarrow X'$. Furthermore, from general
properties of determinantal varieties $Y$ is projectively normal
(see \cite{MS}). We argue as in the proof of \cite[Thm.~2.3]{K} to
show that $\tilde{X}$ is normal. It follows, that the morphism
$Y_P\rightarrow Y$ is given by the linear system $|G|$.
\end{proof}
\begin{rema}
We compute with Singular using the method described in
\cite[Rem.~4.1]{GP} that $h^{1,2}(X)=32$. From the fact that $X'$
has $20$ ordinary double points, we obtain $h^{1,1}(X)=3$.
\end{rema}
In order to see more geometrically the Hodge numbers of $X$ and
$Y$, let us change once more the point of view. Let us first analyze
the K\"{a}hler--Mori cone of a smooth Calabi--Yau threefold $Y$
defined by the $3\times 3$ minors of a $4\times 4$ matrix of
generic linear forms in $\mathbb{P}^{7}$.
 Denote by $S$
 the secant variety of $\mathbb{P}^3 \times \mathbb{P}^3 $
embedded by the Segre embedding in $\mathbb{P}^{15}$. It is well
known (see \cite{Ha}) that the dimension of $S$ is $11$, the
degree $20$ and that $\mathbb{P}^3 \times \mathbb{P}^3 \subset S$
is the singular locus.

\begin{theo}\label{thm 4na4}  The threefold $Z$ in $\mathbb{P}^7$ defined by the
$3\times 3$ minors of a generic $4\times 4$ matrix is a
Calabi--Yau threefold with Picard group of rank $2$. Moreover, the
two faces of the K\"{a}hler--Mori cone give small contractions
into
 nodal complete intersections (with $56$ nodes) of a quadric and a quartic in
$\mathbb{P}^5$.
\end{theo}
\begin{proof} From Proposition \ref{app prop} in the Appendix we obtain $h^{1,1}(Z)=2$.
 The threefold $Z$ can be seen as the intersection of $S$ (the secant variety) with
a generic $7$-dimensional linear subspace $W$ of
 $\mathbb{P}^{15}$. This intersection is a smooth Calabi--Yau
 threefold (this follows from Proposition \ref{st6lemat}).
 The Segre embedding of $\mathbb{P}^3 \times \mathbb{P}^3$ in
$\mathbb{P}^{15}$ is covered by two families of $3$-dimensional
linear spaces (the images of $\{x\} \times \mathbb{P}^3$ and $
\mathbb{P}^3 \times \{x\}$). Let $L_1$ and $L_2$ be two linear
spaces from one of these families and $B$ one from the other
($L_i$
 meets $B$ in one point $l_i$ for $i=1, 2$). The join of
$L_1$ and $L_2$ is a $7$-dimensional linear subspace that contains
all those linear spaces from the family of $L_1$ that correspond
to points on the line $l_1l_2$. We obtain a map $\pi_1\colon
Z\rightarrow G(2,4)$ such that the image of a point $g$ of $Z$ is
the line $l_1l_2$ on $B$, where the linear space spanned by $L_1$
and $L_2$ meets $W$ in $g$. Indeed, since the codimension of
$7$-dimensional spaces that meet $W$ is $4$ (see \cite{Fulton})
the map $\pi_1$ is a birational morphism that contracts a finite
number of disjoint lines (these are found by Schubert cycles
calculations).

We claim that the image of $\pi_1$ is normal. For this it is enough to show
that the singularities appear only when a line is contracted to a
point. This shall follow from an explicit local description
of the map $\pi_1$. Recall first (see \cite[ex.14.16]{Ha}) that
$S\subset \mathbb{P}^{15}$ can be seen as the set of matrices of
rank $\leq 2$ (where $\mathbb{P}^{15}$ is the set of all matrices). A
point $P\in Z$ corresponds to a linear map
$A_P\colon\mathbb{C}^4\rightarrow \mathbb{C}^4$ with
$2$-dimensional image $I_P$ and  $2$-dimensional kernel $K_P$. The
map $\pi_1$ can be seen as a map that associates to a point $P$
the line $I_P\subset \mathbb{P}^3$. Consider the map $$\Omega\colon \mathcal{M}\ni (x^i_j)_{i,j}\longrightarrow [(x^1_1,...,x^1_4),(x^2_1,...,x^2_4)]\in G(2,4),$$ defined on
the subset $\mathcal{M}$ of the set of $4\times 4 $ matrices that
is an appropriate neighborhood of a point $P$ lying outside a
contracted line.

Observe that locally on $Z$ near $P$ the morphisms $\pi_1$ and
$\Omega$ are equal (after an appropriate change of coordinates).
Moreover, if $G$ is a $3$-dimensional linear space, then
$\Omega|_G$ is an isomorphism as soon as it is injective. Thus, it
follows that $\Omega$ restricted to the tangent space $T_P$ to $Z$
at $P$ is an isomorphism (this is the linear subspace of maps in
$W$ carrying the kernel of $A_P$ into the image of $A_P$). Indeed,
since the $7$-dimensional linear space $L_1L_2$ is the set
(outside $\mathbb{P}^3\times\mathbb{P}^3$) of points consisting of
matrices that have common image, and since $W$ meets $L_1L_2$ in
one point, we deduce that $\Omega|_{T_P}$ is injective.

Since $K_Z=0$ and the image of $\pi_1$ is normal, we obtain that
this image is a Calabi--Yau threefold with Gorenstein terminal
singularities (\emph{cDV} singularities). Moreover, since the
Grassmannian $G(2,4)$ is a quadric in $\mathbb{P}^5$ and the image
is normal, we can use Klein theorem \cite[Ex.~6.5]{Har} and
conclude that this image is isomorphic to a complete intersection $X_{2,4} \subset \mathbb{P}^5$
of a quadric and a quartic in $\mathbb{P}^5$.

Let, us show that $X_{2,4}$ is nodal. We claim that the generic
hyperplane section $W$ of $Z$ containing a fixed contracted line
$C$ is smooth. Indeed, observe that it is enough to choose $W$
such that it does not contain any tangent space to $Z$ in points
of the curve $C\subset Z$. Such choice can be done since the
considered tangent spaces to $Z$ induce a curve in the
Grassmannian $G(2,6)$ (of dimension $3$ linear spaces containing
$C$). We can choose for dimension reason $W$ such that the induced
$9$-dimensional $G(2,5)\subset G(2,6)$ (such a family separates points in
$G(2,6)$) is disjoint from the induced curve. We conclude that the
normal bundle of $C\subset Z$ has subbundle $\mathcal{O}(-1)$, and
we can argue as in the proof of \cite[Thm.~2.1]{K}.

To compute the number of nodes, we find the difference between the
Euler characteristics of $Z$ and a smooth complete intersection of
a quadric and a quartic. We obtain that the difference is $112$,
which gives $56$ nodes.
\end{proof}
We next study determinantal varieties using the Grassmann blow up
(see \cite{CM}). Let us describe its properties.

Let $A$ be an $n\times n$ matrix with coordinates as entries. Let
$M_3$ (resp. $M_2$) be the variety given in $\mathbb{P}^{n^2-1}$
by the vanishing of all $3\times 3$ (resp. $2\times 2$) minors of
the matrix $A$. It is a well known fact that $M_2$ is the image of
$\mathbb{P}^{n-1}\times \mathbb{P}^{n-1}$ by its Segre embedding
and that $M_3$ is the secant variety of $M_2$. It follows that
$M_2$ is the singular locus of $M_3$. There are two natural
resolutions of singularities of $M_3$. The first one is the
blowing up of the singular locus $M_{2}$. We will denote this
blowing up by $\pi\colon X\rightarrow M_3$. The exceptional locus
of $\pi$ is a $\mathbb{P}^{n-1}\times \mathbb{P}^{n-1}$ bundle
over $M_2$. The second resolution is given by the formula
\[
 \tau \colon M_3 \times G(n-2,n) \supset Y=\left\lbrace (A,\Lambda) \colon A|_{\Lambda}=0 \right\rbrace \longrightarrow M_3\text{,}
\]
where $\tau$ is the projection onto the first part. The
exceptional locus of $\tau$ is a $\mathbb{P}^{n-1}$ bundle.

\begin{lemm}\label{michal} There is a commutative diagram.
\[\begin{array}{ccc}X&\rightarrow &Y\\
\downarrow& \swarrow&\\
 M_3&&
\end{array}
\]
where the morphism $\sigma$ restricted to the exceptional set over
each point of $M_2$ is the projection onto one of the variables of
$\mathbb{P}^{n-1}\times \mathbb{P}^{n-1}$.
\end{lemm}
\begin{proof}
>From the definition of blowing up, we can view $X$ as a subset of
$\mathbb{P}^{n^2-1}\times
\mathbb{P}^{\left(\frac{n(n-1)}{2}\right)^2-1} $. Using the
identification $\mathbb{P}(M(n,n))=\mathbb{P}^{n^2-1}$, we can
describe $X$ as the set of points of the form $(B,C)$, where $B$
corresponds to a matrix from $M_3$ and $C$ is the
$\frac{n(n-1)}{2} \times \frac{n(n-1)}{2}$ matrix of the $2\times
2$ minors of $B$. The morphism $\pi$ is then the projection onto
the first component of the product.
 Consider now the morphism.
\[
 \tilde{\tau} \colon M_3 \times G(n-2,n) \supset \tilde{Y}=\left\lbrace (A,\Lambda) \colon A^{\textsc{T}}|_{\Lambda}=0 \right\rbrace \longrightarrow M_k\text{.}
\]
We claim that the fiber product $Z=Y\times_{M_3}\tilde{Y}$ is
isomorphic to $X$. In order to prove this, we view $Z$ as a subset
of $$M_3 \times G(n-2,n)\times G(n-2,n)$$ in its turn embedded in
$\mathbb{P}^{n^2-1}\times
\mathbb{P}^{\left(\frac{n(n-1)}{2}\right)^2-1}$ by the composition
of the Pl\"{u}cker embeddings of each Grassmannian and the Segre
embedding. The morphism from $Z$ to $M_3$ is then given by the natural
projection. It is now enough to prove that the traces of $Z$ and
$X$ on the set $M_3\setminus M_2 \times
\mathbb{P}^{\left(\frac{n(n-1)}{2}\right)^2-1}$ coincide, i.e.
that for each matrix $M$ of rank 2 the point
$$\operatorname{Segre}(\operatorname{Pl\ddot{u}cker}(\ker(M)),\operatorname{Pl\ddot{u}cker}(\ker(M^\textsc{T})))$$
has coordinates being the $2\times 2$ minors of $M$. The latter is
checked by direct computation after a suitable choice of
generators of the kernels. This ends the proof of the claim.

It follows that $X$ and $Z$ are the closures of the same set,
hence are equal.\end{proof}

\begin{propo}\label{thm.st 6} The rank of the Picard group of $X$ is 3 and the morphism $\pi \colon X\rightarrow Y$
contracts a two dimensional space of curves on the
K\"{a}hler--Mori cone of $X$. Moreover, the Hodge number $h^{1,2}$
of a generic Calabi--Yau threefold defined by $3\times 3$ minors
of a $4\times 4$ generic matrix is $34$.

\end{propo}
\begin{proof} First by Proposition \ref{app prop} from the Appendix, we obtain that the Chow
group $A^1(S)$ of the secant variety $S$ of $\mathbb{P}^3\times
\mathbb{P}^3\subset \mathbb{P}^{15}$ is
$\mathbb{Z}\oplus\mathbb{Z}$. The blowing up
$\pi\colon\tilde{S}\rightarrow S$ of
$\mathbb{P}^3\times\mathbb{P}^3\subset S$ factorizes through the
Grassmann blow up \cite[p.~206]{Ha} that gives a resolution of
$S$. Denote by $E$ its exceptional set. Since $S$ is regular in
codimension $1$, we obtain from \cite[Prop.~2.6 II]{Har} that the
rank of the Picard group of $\tilde{S}$ is $3$. From the
Grothendieck--Lefschetz theorem applied (several times) to the
pull-back to $\tilde{S}$ of the system of hyperplanes that pass
through a fixed point of $\mathbb{P}^3\times \mathbb{P}^3\subset
\mathbb{P}^{15}$ (this system is very ample), we obtain that the
Picard group of $X$ has rank $3$.

In order to compute the Hodge numbers of the generic element of
the smoothing family of $Y$, we shall use
\cite[Thm.~10]{namikawa1}. We need a description of the
K\"{a}hler--Mori cone of $X$. The hyperplane $W$ defining $Y$ passes
through exactly one point $Q$ from $\mathbb{P}^3 \times
\mathbb{P}^3$ ($Q$ corresponds to the point $(1,0,\dotsc,0)$ in
$Y\subset \mathbb{P}^7$). The incidence correspondence
$$\mathcal{C}=\{(A, \Lambda)\colon A|_{\Lambda}=0 \}\subset Y\times
G(2,4)$$ (where a point of $Y$ corresponds to a matrix $A$ that
gives a linear map $\mathbb{C}^4\rightarrow \mathbb{C}^4$) gives a
partial resolution $\rho_1\colon\mathcal{C}\rightarrow Y$ of $Y$
such that $\mathcal{C}$ is normal (since its image $Y$ is normal).
Observe that the exceptional locus of this Grassmann resolution is
isomorphic to $\mathbb{P}^2$. From Lemma \ref{michal} the blowing
up $X\rightarrow Y$ factorizes through $\rho_1$. We hence obtain a
morphism $\theta_1\colon X\rightarrow \mathcal{C}$ that maps the
exceptional del Pezzo surface $D$ to $\mathbb{P}^2$. It follows
that $\theta_i|_{D}$ for $i=1, 2$ is a blowing down of three
rational curves. Since $K_{X}=0$ these curves map to three
terminal singularities. Since they are contained in a smooth
surface $D$, their normal bundle is $\mathcal{O}(-1)\otimes
\mathcal{O}(-1)$ (see the proof of \cite[Thm.~2.1]{K})). It
follows that these singularities are three ordinary double points
(see \cite[Rem.~5.13(b)]{Reid 3}), thus $\mathcal{C}$ is a nodal
Calabi--Yau threefold.

We see that the image of the restriction map
$$Pic(X)\otimes\mathbb{C}\rightarrow Pic(D)\otimes\mathbb{C}$$ has
dimension $2$. Moreover, the \emph{Case 4} with $r=2$ in Theorem
$10$ from \cite{namikawa1} holds. We deduce that the image of the
natural map of Kuranishi spaces $Def(Y)\rightarrow Def(Y,P)$
coincides with the two dimensional smoothing component $S_1$ of
$Def(Y,P)$. Let $D_{1,loc}$ be the sub-functor of $Def(Y,P)$
corresponding to $S_1$ (see \cite[Lemm.~11]{namikawa1}), then we
have a surjection on tangent spaces $T_{Def(T)}\rightarrow
T_{D_1}$ and we can argue using \cite[Thm.~1.9]{G1} (as in the
proof of Theorem 3.3 \cite{GM}). We obtain that the Hodge number
$h^{1,2}$ of a generic Calabi--Yau threefold defined by $3\times
3$ minors of a $4\times 4$ generic matrix is $$
h^{1,2}(X)+dim(S_1)=32+2.$$ Here we use the fact that $X'$ has
$20$ nodes and that a generic Calabi--Yau threefold defined by
Pfaffians of a $7\times 7$ matrix have $h^{1,1}=1$ and $h^{1,2}
=50$ (see \cite{Rod}).\end{proof}
\begin{rema} From the proof of Lemma \ref{michal}, we obtain the following commutative diagram

             \[ \begin{CD}
             {X}   @>{\theta_2}>>    {\mathcal{D} }\\
             @V{\theta_1}VV         @V\rho_2VV\\
             {\mathcal{C}}   @>\rho_1>> {Y}
             \end{CD}\]
where
$$\mathcal{D}=\{(A, \Lambda)\colon A^r|_{\Lambda }=0 \}\subset Y\times
G(2,4),$$$\theta_1$ and $\theta_2$ are primitive contraction of type
III onto nodal Calabi--Yau threefolds $\mathcal{C}$ and
$\mathcal{D}$. Observe that $\theta_1 \circ \rho_1$ is then the ordinary blow up.
\end{rema}
\section{Del Pezzo of degree $7$} \label{degree7} Let $D'\subset \mathbb{P}^7$ be an anti-canonically
embedded del Pezzo surface of degree $7$. It is well known that
$D'$ can be described by the $2\times 2$ minors of the $3\times 4$
matrix obtained by deleting the last row from a symmetric $4\times
4$ matrix $M$.

We embed $D'$ into a Calabi--Yau threefold $X'$ defined by the
$3\times 3$ minors of the $4\times 4$ matrix

            $$\left( \begin{array}{c}
            \begin{array}{cccc}
            s_1&s_2&s_3&s_4
            \end{array}\\\\

            \begin{array}{c}
            \text{\Huge{\emph{M}}}
            \end{array}
            \end{array}\right)$$

where $s_1,\dotsc,s_4$ are generic linear forms on $\mathbb{P}^7$.

To prove that $X'$ is a nodal Calabi--Yau threefold and to compute
the number of nodes consider the following $4 \times 5$ matrix $N$
  $$ \left( \begin{array}{c|cccc}
\begin{array}{c}
l\\
\end{array}
&
\begin{array}{cccc}
s_1&s_2&s_3&s_4\\
\end{array}\\
\hline
\begin{array}{c}
s_1\\
s_2\\
s_3\\
\end{array}

& \text{\Huge{\emph{M}}}

\end{array} \right)$$

where $l$ is a generic linear forms in $\mathbb{P}^7$ chosen such
that this matrix is obtained from a symmetric $5\times 5$ matrix
by deleting the last row. The $3\times 3$ minors of the matrix $K$
define a smooth canonically embedded surface $G'\subset X'$ of
degree $27$. We can compute using Singular that $G'$ and $D'$ have
$11$ points in common. This suggests the following theorem.
\begin{propo}\label{nody5na4} The threefold $X'\subset \mathbb{P}^7$ is a nodal Calabi--Yau threefold
with $11$ nodes. The blowing up of $D'\subset X'$ is a small
resolution. Let $X\rightarrow X'$ be the flopping of the exceptional
curves of this resolution. Then the Calabi--Yau threefold $X$ contains a
del Pezzo surface $D\simeq D'$ and has Picard group of rank $3$.
\end{propo}
\begin{proof}We argue as in the proof of Theorem \ref{nodal singular}.
Let $Y\subset \mathbb{P}^8$ be the variety defined by the
vanishing of the $3\times 3$ minors of the matrix $\tilde{N}$ obtained
by replacing the upper left entry $l$ in the matrix $N$ by the remaining free
variable. We prove first that $Y$ has exactly one singular point
at $(1,0,\dotsc,0)$, which is resolved by the blowing up
$\mathbb{P}^7\times \mathbb{P}^8 \supset Z\rightarrow Y$ of this point. The
exceptional divisor of this resolution is isomorphic to $D'$. Next, we shall see
that the projection $Z\rightarrow X'\subset \mathbb{P}^7$ is a
small resolution such that $Z=X$.

Let $T$ be the determinantal variety in
$\mathbb{P}^{13}$ of $4 \times 5$ matrices of rank $\leq 2$. The
incidence variety
$$\mathcal{E}=\{(A, \Lambda)\colon A|_{\Lambda}=0 \}\subset T\times
G(2,4)$$ gives a partial resolution $\mathcal{E}\rightarrow T$.
Let us show that $\mathcal{E}$ is smooth. By straightforward
computations, we see that the fibres of the projection $\xi\colon
\mathcal{E}\rightarrow G(2,4)$ are $4$-dimensional projective
spaces in $\mathbb{P}^{13}$. Since $\xi$ is flat, we obtain a
morphism
$$\chi\colon G(2,4)\rightarrow G(4,14).$$ We shall show that the
image of this morphism is smooth. First it is clearly injective,
so the image is generically smooth. Now if $\Lambda_1,\Lambda_2 $
are linear surfaces in $\mathbb{C}^4$ containing $0$, then we can
find an automorphism $P\in GL(4)$ such that
$P\Lambda_1=\Lambda_2$. We obtain an automorphism $\sigma$ of
$G(4,14)$ induced by the linear map $A\rightarrow QAP^{-1}$
between $4\times 5$ partially symmetric matrices, where $Q$ is a
$5\times 5$ matrix with $(P^{-1})^t$ in the upper left corner, $1$
in the lower right corner, and $0$ elsewhere. The automorphism
$\sigma$ maps $\chi(\Lambda_1)$ into $\chi(\Lambda_2)$ and preserves
the image of $\chi$. It follows that the image of $\chi$ is smooth
thus we can argue as in \cite[p.~205]{Ha} and show that
$\mathcal{E}$ is smooth. We conclude that the blowing up $S
\rightarrow T$, that factorizes through $\mathcal{E}\rightarrow T$
(see Lemma \ref{michal}), gives a resolution of $T$. It follows
that the blowing up $Z\rightarrow Y$ in $(1,0,\dotsc,0)$ is a
resolution.

To show that $X'$ is nodal, $Z\rightarrow X'$ is a small
resolution, and that $Z=X$, we argue as in Theorem \ref{nodal
singular}, using the fact that $T$ has Gorenstein singularities
(see \cite{Co}). Finally, from Propositon \ref{app prop} in the
Appendix, we deduce as in the proof of Proposition \ref{thm.st 6}
that $\rho(X)=3$.
\end{proof}
Denote by $G$ the strict transform of $G'$ on $X$.
\begin{propo}\label{lem st7}  The linear system $|G|$ gives a birational morphism
$\pi\colon X\rightarrow Y$ into a
normal variety in $\mathbb{P}^8$ described by the vanishing of the $3\times 3$ minors of
a $4\times 5$ partially symmetric matrix $\tilde{N}$ (from the
proof of Proposition \ref{nody5na4}).The exceptional locus of this morphism is $D$. Moreover, $\pi$ factorizes as
$X\xrightarrow{\rho}V\xrightarrow{\psi} Y$, where $\rho$ is a small
contraction from $X$ into a nodal Calabi--Yau threefold with two nodes.
\end{propo}
\begin{proof}  Let $H$ be the pull back to $X$ of the system of hyperplanes in $\mathbb{P}^8$. We
claim that $G\in |H+D|$. Indeed, let $q$ be the determinant of a
$2\times 2$ minor $B$ of $M$. Then the quadric $q=0$ cuts $X'$ along the
divisor $D'+S'$. Now, by applying the algorithm computing the quotient of ideals,
using an algebra computer system, we show that the $3 \times 3$ minor
of the matrix $K_l$ containing $B'$ with first row and column of
$K_l$ added, determines a cubic that cuts $X$ along $S+G$.

The tangent cone of $Y$ in the point $(1,0,\dotsc,0)\in Y$ is determined by the
vanishing of $2\times 2$ minors of the matrix obtained from
$\tilde{N}$ by deleting the first row and column. Hence, the
exceptional divisor of the blowing up $Z\rightarrow Y$ is
isomorphic to the del Pezzo surface $D'$ (see \cite[p.~257]{Ha}).

We obtain the factorization of $\pi$ from the composed morphism
$S\rightarrow \mathcal{E}\rightarrow T$ (see the proof of Proposition \ref{nody5na4}). The exceptional divisor
of $\psi\colon V\rightarrow Y$ is thus isomorphic to
$\mathbb{P}^2$. We conclude that $V$ has two nodes (see the proof
of Proposition \ref{thm.st 6}).
\end{proof}
\begin{theo} The threefold $R$ defined in $\mathbb{P}^8$ by a generic partially symmetric $4\times 5$ matrix is a smooth
Calabi--Yau threefold with Picard group of rank $2$. Moreover, one
face of its K\"{a}hler--Mori cone gives a contraction with
exceptional set being the surface $\mathbb{P}^1\times
\mathbb{P}^1$ (the image is described in Proposition
\ref{5na5prop}). The other face gives a small contraction to a
nodal Calabi--Yau threefold with $63$ nodes, that is a complete
intersection of a quadric and a quartic in $\mathbb{P}^5$.
Moreover, the Hodge number $h^{1,2}(R)=25$.
\end{theo}
\begin{proof} From the proof of Proposition \ref{nody5na4} the threefold $K$
is smooth (the codimension of the singular locus of $T$ is $4$).
>From Proposition \ref{app prop} in the Appendix, we obtain
$A^1(T)=\mathbb{Z}\oplus \mathbb{Z}$, so from the
Grothendieck--Lefschetz theorem $\rho(R)=2$.

We compute $h^{1,2}(R)$ using the morphism $\pi\colon X\rightarrow
Y$. Indeed, from the Proposition \ref{nody5na4} we find that the
image of the restriction map $Pic(X)\rightarrow Pic(D)$ is
generated by $K_D$ and $E_1+E_2$ (where $E_1$ and $E_2$ are
exceptional divisors on $D$). We obtain, as in
\cite[Thm.~10]{namikawa1}, that $Def(Y)\rightarrow Def(Y,P)$ has
image being the one dimensional smoothing component of
($Def(Y,P)\simeq \mathbb{C}[[x_1,x_2]]/(x_1^2,x_1x_2)$). Now,
arguing as in the proof of Proposition \ref{thm.st 6}, we obtain
$h^{1,2}=25$.

Set $\mathcal{C}=\{(A, \Lambda)\colon A|_{\Lambda}=0 \}\subset
R\times G(2,4)$. Since $R$ is smooth, the natural projection
$\mathcal{C}\rightarrow R$ is an isomorphism. Consider the second
natural projection $p_1\colon \mathcal{C}\rightarrow G(2,4)$ and
the projection $p_2\colon \mathcal{E}\rightarrow G(2,4)$ from the
proof of Proposition \ref{nody5na4}. Observe that $p_1$ is
obtained from $p_2$ by choosing a general linear space
$\mathbb{P}^7\subset\mathbb{P}^{13}\supset \mathcal{E}$. Note that
the fibers of $p_2$ are $4$-dimensional linear space in
$\mathbb{P}^{13}$, forming a $4$-dimensional family
$\chi(G(2,4))\subset G(4,14)$. It follows that $p_1$ is birational
into its image and has a finite number of fibers being lines in
$\mathbb{P}^{7}$. We claim that these lines are contracted to ODP
singularities. To see this, we argue as in the proof of Theorem
\ref{thm 4na4}, showing that we can find a smooth surface in $R$
containing a contracted line. We thus conclude that the image of
$p_1$ is a nodal intersection of a quadric and a quartic in
$\mathbb{P}^5$. To find the number of nodes, we compute the
difference between the Euler characteristics of $R$ and a smooth
complete intersection of a quadric and a quartic in
$\mathbb{P}^5$. The contraction corresponding to the second ray is
discussed in Proposition \ref{5na5prop}.
\end{proof}
\begin{rema}
 Observe that there is another del Pezzo surface of degree $7$ contained in $X'$.
Indeed, consider the surface $\tilde{D}$ defined by the $2\times
2$ minors of the $3\times 4$ matrix obtained from $N$ by deleting
the first and the last row. Now, we can repeat the above
construction for $\tilde{D}\subset X'$, and obtain another
(birational to $Y$) singular Calabi--Yau threefold $Y'$. This one being determined
by the vanishing of the $3\times 3$ minors of a partially
symmetric $5\times 4$ matrix.
\end{rema}

\section{Del Pezzo of degree 8} \label{deg8}
 Let $D\subset \mathbb{P}^8$ be an anti-canonically embedded del
Pezzo surface of degree $8$. We have two possibilities. The surface $D$ is
isomorphic either to $\mathbb{P}^1\times \mathbb{P}^1$ or to
$\mathbb{P}^2$ blown up in one point.
\subsection{}
If the surface $D$ is isomorphic to $\mathbb{P}^1\times
\mathbb{P}^1$, it can be described by the vanishing of the
$2\times 2$ minors of a symmetric $4\times 4$ matrix $M$. We can
as before embed $D$ into a smooth Calabi--Yau threefold $X$
defined by the vanishing of the $3\times 3$ minors of a $5\times
4$ matrix $N$. Indeed, we construct $N$ from $M$ by adding one row
with generic linear forms in such a way that the matrix $N$ could also be
obtained by deleting the first column from a symmetric matrix
$K_l$ (where $l$ is a generic linear form in the upper left corner
of $K_l$). The $3\times 3$ minors of $K_l$ define then a smooth
(because $K_l$ is generic) surface $G$ (we obtain in fact a
family). Denote by $K$ the $5\times 5$ symmetric matrix obtained
from $K_l$ by replacing $l$ by the remaining free variable in
$\mathbb{P}^9$ (as in the proof of Theorem \ref{nodal singular}).
\begin{propo}\label{5na5prop}  The linear system $|G|$ defines a birational morphism with $D$ as exceptional locus, into
a normal Calabi--Yau threefold $Y$ in $\mathbb{P}^9$ defined by
the vanishing of the  $3\times 3$ minors of the symmetric $5\times
5$ matrix $K$ (defined above).
\end{propo}
\begin{proof} First, we see as in the proof of Proposition \ref{lem st7} that $G\in |H+D|$.
Next, we argue as in Proposition \ref{nody5na4}. Let $T\subset
\mathbb{P}^{14}$ be the variety parametrizing symmetric $5\times 5$ matrices of
rank $\leq 2$. Consider the variety
$$\mathcal{E}=\{(A, \Lambda)\colon A|_{\Lambda}=0 \}\subset T
\times G(3,5).$$ Observe that $\mathcal{E}$ is smooth. Now, the
blow up of $T$ along the locus of matrices of rank $1$ has
$\mathbb{P}^3\subset \mathbb{P}^{9}$ (defined by the $2\times 2$
minors of a symmetric $4\times 4$ matrix) as exceptional fibers,
whereas the projection $\alpha \colon \mathcal{E}\rightarrow T$ has
$G(3,4)\simeq \mathbb{P}^3$ as exceptional fibers. So arguing as
in Lemma \ref{michal}, we obtain that $\mathcal{E}\rightarrow T$
is the blowing up of the singular locus of $T$.

Observe that $Y$ is obtained by taking a general linear section of
dimension $9$ passing through a singular point of $T\subset
\mathbb{P}^{14}$. Since $X$ is smooth, the projection
$$\mathcal{E}\supset \alpha^{-1}(Y)\rightarrow G(2,5)\subset \mathbb{P}^9$$ is an
isomorphism onto $X$. Thus the blow up $\alpha^{-1}(Y)\rightarrow Y$
is the morphism given by $|G|$.
\end{proof}
\begin{theo}\label{image} The threefold defined in $\mathbb{P}^9$ by the vanishing of $3\times 3$ minors of a generic symmetric $5\times 5$ matrix is a smooth
Calabi--Yau threefold with Picard group of rank $1$ (i.e.
$h^{1,1}=1$). Moreover, its Hodge number $h^{1,2}=26$.
\end{theo}
\begin{proof}
First,  $X$ can be seen as a linear section of the natural
embedding in $\mathbb{P}^{14}$ of the quotient of
$\mathbb{P}^4\times \mathbb{P}^4$ by the involution
$(x,y)\rightarrow (y,x)$. Since the Picard group of
$\mathbb{P}^4\times\mathbb{P}^4$ is $\mathbb{Z}\oplus \mathbb{Z}$
and the involution transforms one generator into the other, we obtain
that the Picard group of the quotient is $1$. Now, from the
Grothendieck--Lefshetz theorem, we deduce that $\rho(X)=1$. The above
fact is also proved in a more general context in Proposition
\ref{app prop} from the Appendix. To compute the Hodge number $h^{1,2}$, we argue as
in Proposition \ref{thm.st 6}.
\end{proof}
\begin{rema} We shall describe a natural relation between the threefold $X$
and a quintic in $\mathbb{P}^4$, that closes our cascade. As it was
observed in the proof of \cite[Thm.~7.4]{GP}, the smooth
Calabi--Yau threefold $X$ defined by the $3\times 3$ minors of a
symmetric $5\times 5$ matrix in $\mathbb{P}^9$ admits an unramified
covering being a Calabi--Yau threefold. Indeed, let $T$ be the
pre-image to $\mathbb{P}^4\times \mathbb{P}^4$ of $X$ by the
involution.
 The image of the projection $p_1$ of $T$ on $\mathbb{P}^4$
can be seen to be a quintic (the determinant of a $5\times 5$
matrix) with $50$ nodes (for generic choice of $X$). The
projection $T\rightarrow p_1(T)$ is then a primitive contraction of
type I.

\end{rema}
\subsection{}
 If the surface $D$ is isomorphic to $\mathbb{F}_1$ ($\mathbb{P}^2$ blown up in one
point), it can be described by the $2\times 2$ minors of a
$3\times 5$ matrix $R$ of linear forms in $\mathbb{P}^8$, obtained
by joining two symmetric $3\times 3$ matrix with one common
column.
$$\left( \begin{array}{ccc}
\begin{array}{cc|c|cc}
l_1&l_2&l_3&l_4&l_5\\
l_2&l_6&l_4&l_7&l_8\\
l_3&l_4&l_5&l_8&l_9\\
\end{array}

\end{array}\right)$$

We will call such matrices double-symmetric. Let us embed $D$ into a
singular Calabi--Yau threefold $X'$ defined by the $3\times 3$
minors of the $4\times 5$ matrix $T$ obtained by adding one row to
the matrix $R$ in such a way that $T$ can be obtained by deleting
the last row from a symmetric $5\times 5$ matrix. Denote by
$Q_{l_{12}}$ (or $Q$) the following double-symmetric $4\times 6$
matrix
$$\left( \begin{array}{ccc}
\begin{array}{cc|cc|cc}
l_1&l_2&l_3&l_4&l_5&l_8\\
l_2&l_6&l_4&l_7&l_8&l_{10}\\
l_3&l_4&l_5&l_8&l_9&l_{11}\\
l_4&l_7&l_8&l_{10}&l_{11}&l_{12}
\end{array}
\end{array}\right)$$
where $l_{10},l_{11},l_{12}$ are generic linear forms in
$\mathbb{P}^8$. Let $G_{l}\subset X'$ be the surface defined by
the $2\times 2$ minors of $Q_l$.
\begin{propo} The threefold $Y$ defined by the $3\times 3$ minors of a
generic $4\times 6$ double-symmetric matrix of linear forms in
$\mathbb{P}^9$ has $12$ isolated singular points analytically
isomorphic to cones over $\mathbb{F}_1$. Moreover, the threefold
$X'$ has $1$ ordinary double point and $11$ more singularities
described below.
\end{propo}
\begin{proof} Consider the variety $T_k\subset \mathbb{P}^{11}$ parametrizing double symmetric
$4\times 6$ matrices of rank $\leq k$. Choosing coordinates
$l_1,\dotsc,l_{12}$ in $\mathbb{P}^{11}$, the scheme $T_k$ is
defined by the $k+1\times k+1$ minors of $Q$. We find using
Singular that the dimensions of $T_2$ and $T_1$ are $5$ and $2$
(and degrees $12$ and $35$) respectively. Let us show that
$T_2-T_1$ is smooth and that there exists a linear change of
coordinates of $\mathbb{P}^{11}$ preserving $4\times 6$ double
symmetric matrices that transforms a given element of $T_1$ to
$(1,0,\dotsc,0)$. In particular the type of singularities on $T_2$
wll then be the same in all points of $T_1$.
 Observe first that the matrix $Q$ can
be described in the following form
$$\left( \begin{array}{ccc}
\begin{array}{ccc}
A&B&C\\
B&C&D\\
\end{array}
\end{array}\right)$$
where $A$, $B$, $C$, and $D$ are symmetric $2\times 2$ matrices.
Consider the following operations preserving double symmetric
matrices.
\begin{enumerate}  \item
The transformation that
change $(A,B,C,D)$ from the matrix $Q$ into
$$ (A,sA+B,(s^2+s)A+sB+sC,2(s^3+s^2)A+(2s^2+s)B+2sC),$$ for chosen $s\in \mathbb{C}$.
\item The central symmetry. \item The operations between symmetric $2\times 2$ matrices
$$\left(
\begin{array}{cc}
a,b\\
b,c\\
\end{array}\right)
\longrightarrow \left(\begin{array}{cc}
a,rb+ta\\
b+sa,rc+sb+s^2a\\
\end{array}
\right), $$ for chosen $t,r\in \mathbb{C}$, performed
simultaneously on $A,B,C,D$.
 \item The operation of exchanging
rows $1$ with $2$ and $3$ with $4$, composed with the operation exchanging columns $1$ with
$2$, $3$ with $4$, and $5$ with $6$.\end{enumerate}

We claim that the composition group of the above operations acts transitively on $T_1$. Indeed, the
rank $1$ double symmetric matrices are exactly those with $aA=b
B=cC=dD$ such that $a,b,c,d\in\mathbb{C}$,
$ac-b^2=bd-c^2=ad-bc=0$, and such that $A,B,C,D$ have rank $1$.

Observe that such $Q\in T_2-T_1$ with $B=0$ and $aA=cC=dD$ of rank $1$
have two dimensional orbits of the action of the group. Except in the latter case, we can
find operations that transform a matrix $Q\in T_2-T_1$ into a
matrix $R$ with $A$ and $B$ of rank $2$. Since $R$ has rank $2$,
we compute that $AC-B^2=AD-BC=BD-C^2=0$. In fact if such equation
is satisfied and at least one of the matrices $A,\ B,\ C,\ D$ have
rank $2$ then $R$ has rank $2$. In the above case the orbits of
the considered operations are three dimensional.

In a neighborhood of $R\in T_2-T_1$ with $A$ and $B$ of rank $2$,
we find a natural parametrization (fixing the $6$ entries of $A$ and $B$)
and conclude that $T_2$ is smooth in $R$. To show the smoothness
in the points corresponding to matrices $Q\in T_2-T_1$ with $$A= \left( \begin{array}{cc}
1&0\\
0&0\\
\end{array}
\right),\ A=aC,\ D=dC $$, where $a\neq 0$ (the case $C=B=0$ is
analogical), consider the intersection $T_2\cap H_2\cap H_3\cap
H$, where $H_i=\{ l_i=0 \}$
 for $i=2,3$ and $H=\{ l_5=d\cdot l_9 \}$. We find a local parametrization of this intersection in a neighborhood of $Q$
 with complex plane coordinates $(x,r)$ close to $(a,0)$
$$\left( \begin{array}{c}
\begin{array}{cccccc}
x&0&0&xr&1&dxr\\
0&x^2r^2&xr&dx^2r^2&dxr&d^2x^2r^2+xr^2\\
0&xr&1&dxr&d&d^2xr+r\\
xr&dx^2r^2&dxr&d^2x^2r^2+xr^2&d^2xr+r&d^3x^2r^2+2dxr^2\\
\end{array}
\end{array}\right). $$
Since $T_2\cap H_2\cap H_3\cap H$ is smooth, we obtain that
$T_2-T_1$ is smooth.

We claim that the blow up of $T_1\subset T_2$ gives a resolution
of $T_2$. Indeed, we know from the descriptions above that $T_1$
is smooth and the generic three dimensional and transversal to
$T_1$ complete intersection in $T_2$ has a node as singularity.

We obtain that the singularities of $Y$ are locally isomorphic to a cone
over $\mathbb{F}_1$ (the $9$-dimensional linear section of $T_2$
in $(1,0,\dotsc,0)$ has tangent cone being a cone over
$\mathbb{F}_1$).

Using Singular, we compute that for generic $l$ the multiplicity
of the intersection of surfaces $G_l$ and that $D$ is $34$ and the
radical has degree $12$. This suggests the following description.
Consider the threefold $Y\subset \mathbb{P}^9$ defined by $3\times
3$ minors of $Q_x$, where $x$ in a new free variable. Blow up
$\mathbb{P}^8\times\mathbb{P}^9\supset Z\rightarrow Y\subset
\mathbb{P}^9$ in $(0,\dotsc,0,1)$. The projection $Z\rightarrow
\mathbb{P}^8$ contracts $12$ lines $11$ of which pass through the
singular points of $Z$. We see that the remaining line has normal
bundle $\mathcal{O}(-1)\oplus \mathcal{O}(-1)$ and is contracted
to an ODP.
\end{proof}
\begin{rema} Observe that if we continue the cascade and embed
the del Pezzo surface of degree $9$, that is defined by
$2\times2\times 2$ minors of a $3\times3\times 3$ matrix with
linear forms in $\mathbb{P}^9$ symmetric with respect to three
rectangular diagonals containing a chosen main diagonal, into one
of the above Calabi--Yau threefolds, we obtain a
threefold that cannot be smoothed (his singularity is rigid). The resulting variety
is possibly defined by $3\times 3\times 3$ minors of a
$4\times4\times 4$ matrix symmetric with respect to three
rectangular diagonals containing a chosen main diagonal.
\end{rema}

\vskip10pt
\begin{minipage}{4.5cm} Grzegorz Kapustka \\
Jagiellonian University\\
ul. Reymonta 4\\
30-059 Krak\'{o}w\\
Poland\\

\end{minipage}
\hfill
\begin{minipage}{4.5cm}
 Micha\l\ Kapustka\\
 Jagiellonian University\\
 ul. Reymonta 4\\
 30-059  Krak\'{o}w\\
 Poland\\

\end{minipage}\\
email:\\
Grzegorz.Kapustka@im.uj.edu.pl\\
Michal.Kapustka@im.uj.edu.pl

\newpage
\pagestyle{plain}
\section*{Appendix}
\bigskip
\begin{center}\textbf{A note on the Chow groups of projective determinantal varieties}

\end{center}

\bigskip
\begin{center}{by Piotr Pragacz}
\end{center}
 \bigskip

In the present note we shall consider the following types of determinantal varieties.

\noindent
(i) (generic) Let $W{,}V$ be two vector spaces over an arbitrary field $K$ with
$m=\mbox{dim} \ W \ge n= \mbox{dim} \ V$.
For $r\ge 0$, set
\begin{equation}\label{Dr}
D_r = D_r(\varphi)=\{x\in P: \mbox{rank} \ \varphi(x)\le r\}\,,
\end{equation}
where $\varphi: W_P \longrightarrow V_P \otimes \mathcal{O}(1)$
is the canonical morphism on $P=\mathbb{P}(\mbox{Hom}(W{,}V))$.

\noindent
(ii) (symmetric) Take the following specialization of (i): let $m=n$, $W=V^*$, $P=\mathbb{P}(\mbox{Sym}^2(V))$,
and $\varphi: V^*_P \longrightarrow V_P \otimes \mathcal{O}(1)$ be the canonical symmetric
morphism on $P$. Define $D_r$ by (\ref{Dr}).

\noindent
(iii) (partially symmetric) Consider the following specialization of (i): let $m>n$,
$W^*\twoheadrightarrow V$, $P=\mathbb{P}(W^*\vee V)$ (in the notation of \cite{LP}),
and let
$$
\varphi: W_P \longrightarrow V_P \otimes \mathcal{O}(1)
$$
be the canonical partially symmetric morphism on $P$. Define $D_r$ by (\ref{Dr}).

\noindent
(One can also, for even $r$, consider the skew-symmetric analogs of (ii) and (iii).)

\noindent
In all cases (i), (ii) and (iii), we get a sequence of {\it projective determinantal varieties}
$$
\emptyset = D_0 \subset D_1 \subset D_2 \subset \cdots \subset D_{n-1} \subset D_{n}=D_{n+1}=\cdots .
$$
The scheme $D_r$ can be seen as a variety defined in $P$
by the vanishing of $(r+1)\times (r+1)$ minors of a generic $m\times n$
matrix of linear forms. The codimensions of the determinantal varieties $D_r$ in the
respective cases are:
(i) \ $(m-r)(n-r)$; (ii) \ $(n-r)(n-r+1)/2$; (iii) \ $(m-n)(n-r) + (n-r)(n-r+1)/2$.

In the present note, we compute $A^1(D_r)$ for the above determinantal varieties,
getting the answers: $\mathbb{Z}\oplus \mathbb{Z}$ in cases (i) and (iii), and\ $\mathbb{Z}$
in case (ii).
We also discuss generators of the Chow groups $A_*(D_r\setminus D_{r-1})$ in case (i);
for $r=1$ and $r=n-1$, we give some linearly independent generators.

\textbf{Background.} The content of this note was obtained in the late 80's,
and has not been written up to now. Due to a recent ask of G. and M.~Kapustka,
we have decided to write this material up because it is needed in their research.

In this note, we shall use notation, conventions, and some results
from \cite{P}, \cite {LP}. In particular, as for what concerns the
Chow groups, we shall use notation and conventions from \cite{F}.
\setcounter{section}{0}
\section{$D_r\setminus D_{r-1}$ as fiber bundles}

In this section, we follow basically \cite{P}.

\noindent
(i) \ For $f\in \mbox{Hom}(W{,}V)$, we
set $K_f=\mbox{Ker}(f)$, $C_f= \mbox{Coker}(f)$. When $f$ varies in $D_r\setminus D_{r-1}$,
we get the vector bundles $K$ and $C$ of ranks $m-r$ and $n-r$ on $D_r\setminus D_{r-1}$.
We consider the fibration
$$
D_r\setminus D_{r-1}\longrightarrow G_{m-r}(W)\times G_r(V)
$$
given by $f \mapsto (K_f,\ C_f)$. Its fiber is equal to the space of nonsingular
$r\times r$ matrices over $K$. More explicitly, let
$$
P'= \mathbb{P}(\mbox{Hom}(Q_W{,}R_V))\longrightarrow G_{m-r}(W)\times G_r(V)\,.
$$
The bundle $Q_W$ is the pullback on $G_{m-r}(W)\times G_r(V)$ of the
tautological quotient rank $r$ bundle on $G_{m-r}(W)$. Moreover, the bundle $R_V$ is the
pullback on $G_{m-r}(W)\times G_r(V)$ of the tautological subbundle on $G_r(V)$.

On $P'$, there is the tautological morphism
$$
\varphi': (Q_W)_{P'}\longrightarrow(R_V)_{P'}\otimes\mathcal{O}_{P'}(1)\,,
$$
and we have
\begin{equation}\label{iso}
P'\setminus D_{r-1}(\varphi')\cong D_r\setminus D_{r-1}\,.
\end{equation}

\medskip
\noindent
(ii) \ For symmetric $f\in \mbox{Hom}(V^*,V)$, we have $K_f\cong C_f^*$.
We consider the fibration
$$
D_r\setminus D_{r-1}\longrightarrow G_r(V)
$$
given by \ $f \mapsto C_f$. Its fiber is equal to the space of nonsingular symmetric
$r\times r$ matrices. To be more explicit, let
$$
P'= \mathbb{P}(\mbox{Sym}^2(R))\longrightarrow G_r(V)\,,
$$
where $R$ is the tautological subbundle on $G_r(V)$. On $P'$, there is the tautological
symmetric morphism
$$
\varphi': R_{P'}^*\longrightarrow R_{P'}\otimes\mathcal{O}_{P'}(1)\,,
$$
and we have $P'\setminus D_{r-1}(\varphi')\cong D_r\setminus D_{r-1}$.

\medskip
\noindent
(iii) \ For a partially symmetric $f\in \mbox{Hom}(W,V)$, we have $K_f^*\twoheadrightarrow C_f$.
Let $Fl$ denote the flag variety parametrizing the pairs $(A,B)$, where $A$ is an $(m-r)$-dimensional
quotient of $W^*$, $B$ is an $(n-r)$-dimensional quotient of $V$ and we have $A \twoheadrightarrow B$.
We consider the fibration
$$
D_r\setminus D_{r-1}\longrightarrow Fl
$$
given by $f \mapsto (K_f^*, C_f)$. Its fiber is equal to the space of nonsingular $r\times r$ symmetric
matrices. More explicitly, let
$$
P'= \mathbb{P}(\mbox{Sym}^2(R))\longrightarrow Fl\,.
$$
The bundle $R$ is here the tautological rank $r$ subbundle on $Fl$.
On $P'$, there is the tautological symmetric morphism
$$
\varphi': R_{P'}^*\longrightarrow R_{P'}\otimes\mathcal{O}_{P'}(1)\,,
$$
and we have $P'\setminus D_{r-1}(\varphi')\cong D_r\setminus
D_{r-1}$.
\section{Computations of $A^1(D_r)$}

Let $i'$ denote the embedding $D_{r-1}(\varphi')) \to P'$.

\begin{lemm}\label{lem0} In each case (i), (ii) and (iii), we have the following exact sequence of
the Chow groups:
\begin{equation}\label{seq}
A_*(D_{r-1}(\varphi'))\xrightarrow{i'_*} A_*(P')\longrightarrow A_*(D_r\setminus D_{r-1})\longrightarrow 0\,.
\end{equation}
\end{lemm}
This follows by combining (\ref{iso}) and its analogues with \cite{F}, Sect.1.8 applied to
the embedding $D_{r-1}\subset D_r$.

With the help of the Schur S- and Q-functions (cf, e.g.,
\cite{P}), we now record

\begin{lemm}\label {lem1} In case (i), the image $\mbox{Im}(i'_*)$ is generated by
$$
s_I(Q)-s_I(R\otimes L)\,,
$$
where
$$
Q=(Q_W)_{P'}{,}\quad R=(R_V)_{P'}{,}\quad L=\mathcal{O}_{P'}(1)\,,
$$
and $I$ runs over all partitions of positive weight.

\smallskip

In cases (ii) and (iii), by putting $M$ to be the formal square root of $L$,
the image $\mbox{Im}(i'_*)$ is generated by
$$
Q_I(R\otimes M)\,,
$$
where $R$ denotes the pullback to $P'$ of the corresponding tautological rank $r$
subbundle (on $G_r(V)$ or $Fl$), and and $I$ runs over all (strict) partitions of positive
weight.
\end{lemm}
This follows from \cite{P}, Corollary 3.13 and its symmetric
analog established also in \cite{P}.

\begin{propo}\label{app prop} In cases (i) and (iii), we have $A^1(D_r)\cong \mathbb{Z} \oplus \mathbb{Z}$
for any $r\ge 1$. In case (ii), we have $A^1(D_r)\cong \mathbb{Z}$
for any $r\ge 1$.
\end{propo}
\textbf{Proof.} Since $\mbox{codim}(D_{r-1}, D_r)\ge 2$ for $r\ge 1$, it suffices to prove the same
assertions for $A^1(D_r\setminus D_{r-1})$ instead of $A^1(D_r)$. Set, in all three cases,
$h=c_1(L)$.

\noindent
(i) By Lemmas \ref{lem0} and \ref{lem1} we see that $A^1(D_r \setminus D_{r-1})$ is generated
(over $\mathbb{Z}$) by $s_1(Q)$, $s_1(R)$, and $h$, modulo the following single relation:
$$
s_1(Q)=s_1(R\otimes L)=s_1(R)+h\,.
$$
Thus the assertion follows.

\noindent
(ii) We see that $A^1(D_r \setminus D_{r-1})$ is generated
by $s_1(R)$ and $h$, modulo the following single relation:
$$
Q_1(R\otimes M)=2(s_1(R)+s_1(M))=2s_1(R)+h=0\,,
$$
which implies the assertion.

\noindent
(iii) Since $Fl$ is a Grassmann bundle over a Grassmannian, we have
$$
A^1(Fl)\cong \mathbb{Z} \oplus \mathbb{Z}=\mathbb{Z}s_1(R)\oplus \mathbb{Z}x\,,
$$
for some $x$. We see that $A^1(D_r \setminus D_{r-1})$ is generated
by $s_1(R)$, $x$ and $h$, modulo the following single relation:
$$
Q_1(R\otimes M)=2s_1(R)+h=0\,.
$$
Hence the assertion follows.
\qed

Similarly, one shows that the Chow groups $D_r$ ($r$ even) of skew-symmetric and
partially skew-symmetric projective determinantal varieties are of rank 1 and 2,
respectively.

\section{Remarks on other Chow groups of $D_r\setminus D_{r-1}$}

We work here in the generic case (i).

\begin{propo}\label{prop2} For $r\ge 1$, we have the following inequalities:
\begin{equation}
\binom{n}{r} \le \mbox{rank} \ A_*(D_r\setminus D_{r-1}) \le \binom{n}{r}(m-r+1)\,.
\end{equation}
\end{propo}
\textbf{Proof.} To prove the first inequality, we invoke the
following exact sequence of the Chow groups (cf.~\cite{F}, Example
2.6.2):
\begin{equation}\label{seq1}
A_k(D_r) \xrightarrow{\cdot \ h} A_{k-1}(D_r)
\rightarrow A_k(C_{D_r})\rightarrow 0\,,
\end{equation}
where $C_{D_r}$ is the affine cone over $D_r$. We recall the
following result from \cite{Pap}, Proposition 4.2 (recall that we
assume $m\ge n$):
\begin{equation}\label{ab}
\mbox{rank} \ A_*(C_{D_r})=\binom{n}{r}\,.
\end{equation}
The equality (\ref{ab}), combined with the surjection in the sequence
(\ref{seq1}), implies the first inequality.

To prove the second equality, we show that the elements
$s_I(R)\cdot h^j$, where $I\subset(n-r)^r$ and $j=0, \ldots, m-r$,
generate over $\mathbb{Q}$ the Chow group $A^k(D_r\setminus
D_{r-1})$, where $k=|I|+j$. It follows from Schubert calculus
(cf., e.g.~\cite{F}, Chap.14) and the surjection in (\ref{seq})
that the group $A_*(D_r\setminus D_{r-1})$ is generated by
$s_I(Q)$, $I\subset (r)^{m-r}$; $s_J(R)$, $J\subset (n-r)^r$; and
powers of the class $h$. By Lemma \ref{lem1}, in $A_*(D_r\setminus
D_{r-1})$ we have
$$
s_I(Q)=s_I(R\otimes L)\,,
$$
and we see that the group $A_*(D_r\setminus D_{r-1})$ is generated by $s_J(R)$ (with $J\subset (n-r)^r$)
and powers of the class $h$.

If $I\nsubseteq (n-r)^r$, then $s_I(Q) = 0$. Thus, invoking the Lascoux formula for the Schur
polynomial of the twisted vector bundle (cf., e.g., \cite{F}, Ex. A.9.1), we get for such $I$:
\begin{equation}\label{zero}
0 = \sum_{J\subset I, \ J\subset (n-r)^r} d_{IJ} \cdot s_J(R)\cdot h^{|I|-|J|}\,.
\end{equation}
These relations allow us to express the powers $h^{m-r+1}, h^{m-r+2},\ldots$
with the help of $h^j$, $0\leq j\leq m-r$, and $s_I(R)$, $I\subset (n-r)^r$.
\qed

\begin{ex} Let $r=1$. By the proof of Proposition \ref{prop2}, we know that
$s_i(R)\cdot h^j$, where $i=0,1,\ldots, n-1$ and $j=0,1,\ldots, m-1$, generate
over $\mathbb{Q}$ the Chow group $A_*(D_1)$. But by the Segre embedding, we have
$D_1\cong \mathbb{P}(W) \times \mathbb{P}(V)$. Since
$\mbox{rank} \ A_*(\mathbb{P}(W) \times \mathbb{P}(V))=mn$,
the displayed elements are, in fact, $\mathbb{Z}$-linearly independent
generators of $A_*(D_1)$. This can be also seen from the relations given
in the proof of Proposition \ref{prop2}.
\end{ex}

\begin{ex} Let now $r=n-1$.\footnote{This computation was done in collaboration with S.A.~Str\o mme.}
In this case, $C$ is a line bundle and
$G_r(V)\cong \mathbb{P}^{n-1}$. Set $c=c_1(C)$. From the long exact
sequence of bundles relating $K$ and $C(h)$ one gets, for $a\ge m-n+2$,
\begin{equation}\label{10}
\binom{m}{m-a}\cdot h^a - \binom{m}{m-a+1}\cdot h^{a-1}\cdot c =0\,.
\end{equation}
With the help of (\ref{10}), one can deduce that the elements:
$$
h^i\cdot c^j \ \ (0\le i\le m-n, \ 0\le j \le n-1), \ \ \ \ h^{m-n+1}, \ \ \ \ h^{m-n+1}\cdot c
$$
are $\mathbb{Q}$-linearly independent generators of $A_*(D_{n-1}\setminus D_{n-2})$.
\end{ex}

\bigskip\bigskip\bigskip

\small Institute of Mathematics

\small Polish Academy of Sciences

\small \'Sniadeckich 8

\small 00-956 Warszawa

\small Poland

\bigskip

\small P.Pragacz@impan.gov.pl

\end{document}